\documentclass[12pt]{amsproc}
\usepackage{hyperref}
\title[Divisibility of Character Values]{On the Divisibility of Character Values of the Symmetric Group}
\author[Ganguly]{Jyotirmoy Ganguly}
\address{Indian Institute of Science Educations and Research, Pune.}
\author[Prasad]{Amritanshu Prasad}
\address{The Institute of Mathematical Sciences, Chennai (Homi Bhabha National Institute).}
\author[Spallone]{Steven Spallone}
\address{Indian Institute of Science Educations and Research, Pune.}
\subjclass[2010]{20C30,05A16,05A17}
\keywords{symmetric groups, irreducible characters, divisibility, core towers}
\newtheorem*{theoremm}{Main Theorem}
\newtheorem*{theorema}{Theorem A}
\newtheorem{theorem}{Theorem}
\newcommand{\TT}{\mathcal T}

\begin{document}
\begin{abstract}
  Fix a partition $\mu=(\mu_1,\dotsc,\mu_m)$ of an integer $k$ and positive integer $d$.
  For each $n>k$, let $\chi^\lambda_\mu$ denote the value of the irreducible character of $S_n$ at a permutation with cycle type $(\mu_1,\dotsc,\mu_m,1^{n-k})$.
  We show that the proportion of partitions $\lambda$ of $n$ such that $\chi^\lambda_\mu$ is divisible by $d$ approaches $1$ as $n$ approaches infinity.
\end{abstract}
\maketitle
Let $k$ be a positive integer, and $\mu=(\mu_1,\dotsc,\mu_m)$ a partition of $k$.
For a partition $\lambda$  of an integer $n\geq k$, let $\chi^\lambda_\mu$ denote the value of the irreducible character of $S_n$ corresponding to $\lambda$ at an element with cycle type $(\mu_1,\dotsc,\mu_m,1^{n-k})$.
The purpose of this article is to prove:
\begin{theoremm}
  For any positive integers $k$ and $d$, and any partition $\mu$ of $k$,
  \begin{displaymath}
    \lim_{n\to \infty}\frac{\#\{\lambda\vdash n\mid \chi^\lambda_\mu\text{ is divisible by $d$}\}}{p(n)} = 1.
  \end{displaymath}
  Here $p(n)$ denotes the number of partitions of $n$.
\end{theoremm}
In particular, for any integer $d$, the probability that an irreducible character of $S_n$ has degree divisible by $d$ converges to $1$ as $n\to \infty$.

Recall the theorem of Lassalle \cite[Theorem~6]{MR2368985}, which implies that there exists an integer $A^\lambda_\mu$ such that
\begin{equation}
  \label{eq:lassalle}
  \chi^\lambda_\mu = \frac{f_\lambda}{(n)_k}A^\lambda_\mu.
\end{equation}
Here $(n)_k=n(n-1)\dotsb(n-k+1)$, and $f_\lambda$ is the degree of the irreducible character of $S_n$ corresponding to $\lambda$.
Therefore, in order to prove the main theorem, we focus on the divisibility properties of $f_\lambda$.
For each prime number $q$, let $v_q(m)$ denote the $q$-adic valuation of an integer $m$, in other words, $q^{v_q(m)}$ is the largest power of $q$ that divides $m$.
Also write $\log n=\log_q n$.
The main theorem will follow from the following result:
\begin{theorema}
  For every prime number $q$ and non-negative integer $m$,
  \begin{displaymath}
    \lim_{n\to \infty} \frac{\#\{\lambda\vdash n\mid v_q(f_\lambda)\leq m+(q-1)\log n\}}{p(n)}=0.
  \end{displaymath}
\end{theorema}
In the rest of this article, we first prove Theorem A, and then show that it implies the main theorem.

\section{Proof of Theorem A}
\label{sec:proof-theorem}
The proof of Theorem A is based on the theory of $q$-core towers.  This construction originated in the seminal paper \cite{Macdonald} of Macdonald, and was developed further by Olsson in \cite{Olsson}.  We now recall the relevant aspects.

Let $[q]$ denote the set $\{0,\dotsc,q-1\}$, and
consider the disjoint union
\begin{displaymath}
  T_q = \coprod_{i=0}^\infty [q]^i.
\end{displaymath}
The set $T_q$ can be regarded as a rooted $q$-ary tree with root $\emptyset$.
The children of a vertex $(a_1,\dotsc,a_i)\in [q]^i$ are the vertices $(a_1,\dotsc,a_i,a_{i+1})$, where $a_{i+1}\in [q]$.
A partition $\lambda$ is said to be a $q$-core if no cell in its Young diagram has hook length divisible by $q$.
Denote the set of all $q$-core partitions by $C_q$.
The $q$-core tower construction associates to each partition $\lambda$ of $n$ a function $\TT_q^\lambda:T_q\to C_q$ known as the $q$-core tower of $\lambda$.
For a partition $\lambda$, define:
\begin{displaymath}
  w_i(\lambda) = \sum_{x\in [q]^i} |\TT_q^\lambda(x)|.
\end{displaymath}
Then the $q$-core tower satisfies the following constraint:
\begin{equation}
  \label{eq:1}
  \sum_{i=0}^\infty  w_i(\lambda)q^i = n.
\end{equation}
In particular, $\TT^\lambda_q(x)=\emptyset$ for all $i>\log_q n$.
This function $\lambda\mapsto \TT_q^\lambda$ is a bijection from the set of partitions of $n$ onto the set of $q$-core towers satisfying the condition \eqref{eq:1}.

Let $n$ be a positive integer with $q$-ary expansion:
\begin{equation}
  \tag{$*$}
  \label{eq:q-ary}
  n = a_0 + a_1q + \dotsb + a_rq^r, \text{ with } a_i\in [q] \text{ for } i=1,\dotsc,r, \text{ and } a_r>0.
\end{equation}
Define $a(n)=\sum_{i=0}^r a_i$.

Recall the following Theorem:
\begin{theorem}[{\cite[Equation~(3.3)]{Macdonald}}]
  \label{theorem:olsson}
  For a partition $\lambda$, let $w(\lambda)=\sum_{i=0}^r w_i(\lambda)$.
  For any partition $\lambda$ of $n$ and any prime $q$,
  \begin{displaymath}
    v_q(f_\lambda) = \frac{w(\lambda)-a(n)}{q-1}.
  \end{displaymath}
\end{theorem}
Theorem~\ref{theorem:olsson} can be used to find constraints on partitions with small values of $v_q(f_\lambda)$.
Suppose that $v_q(f_\lambda) \leq b$.
By Theorem~\ref{theorem:olsson}, this is equivalent to
\begin{displaymath}
  w(\lambda)\leq a(n)+b(q-1).
\end{displaymath}
The expansion \eqref{eq:q-ary} implies that $r\leq \log n<r+1$, so that $a(n)\leq (q-1)(r+1)\leq (q-1)(\log n+1)$.
So if $v_q(f_\lambda)\leq b$, then
\begin{displaymath}
  w(\lambda)\leq (q-1)(\log n + 1 + b).
\end{displaymath}
Thus an upper bound for the number $p_b(n)$ of partitions $\lambda$ of $n$ such that $v_q(f_\lambda) \leq b$ can be obtained by counting the number of $q$-core towers with $(q-1)(\log n + 1 + b)$ or fewer cells.
The total number of vertices in the first $r+1$ rows of $T_q$, i.e., in  $\coprod_{i=0}^r [q]^i$, is:
\begin{displaymath}
  1+q+\dotsb+q^r=\frac{q^{r+1}-1}{q-1}<qn,
\end{displaymath}
since $q^r\leq n$.
Let $c_q(n)$ denote the number of $q$-core partitions of $n$.
Set $N_b=(q-1)(\log n+b+1)$.
Let $\tilde c_q(n)$ denote $\max\{c_q(i)\mid 1\leq i\leq n\}$.  There are $\binom{w+N-1}{w}$ ways to distribute $w$ cells into $N$ nodes.  
Thus
\begin{align*}
  p_b(n)&\leq \tilde c_q(N_b)^{N_b}\binom{qn+N_b}{N_b}\\
        &\leq \tilde c_q(N_b)^{N_b} (qn+N_b)^{N_b}\\
\end{align*}
It is known that, for every integer $q$, there exists a polynomial $f_q(n)$ such that $\tilde c_q(n)\leq f_q(n)$ for all $n\geq 0$.
Indeed, for $q=2$, it is well-known that $c_2(n)\leq 1$, and for $q=3$, using a formula of Granville and Ono \cite[Section~3, p.~340]{MR1321575}, $c_3(n)\leq 3n+1$.
For $q\geq 4$, the existence of $f_q(n)$ follows from Anderson~\cite[Corollary 7]{MR2444214}.

We get:
\begin{displaymath}
  p_b(n)\leq f_q(N_b)^{N_b}(qn+N_b)^{N_b},
\end{displaymath}
whence
\begin{displaymath}
  \log p_b(n) \leq  N_b[\log f_q(N_b)+\log(qn + N_b)].
\end{displaymath}
Taking $b=m+(q-1)\log n$ gives $N_b=(q-1)(q\log n + m + 1)$.
Thus $\log p_b(n)=o(n^\epsilon)$ for every $\epsilon>0$.
On the other hand, the Hardy-Ramanujan asymptote \cite{HR} for $p(n)$ implies that $\log p(n)$ grows faster than $n^{\frac 12-\epsilon}$ for any $\epsilon>0$.
Thus Theorem~A follows.

\section{Proof of the Main Theorem}
\label{sec:proof-main-theorem}
The identity \eqref{eq:lassalle} implies that
\begin{displaymath}
  v_q(\chi^\lambda_\mu)\geq v_q(f_\lambda)-v_q((n)_k).
\end{displaymath}
Using Legendre's formula on the valuation of a factorial, that $v_q(n!)=\frac{n-a(n)}{q-1}$, we have:
\begin{displaymath}
  v_q((n)_k) = v_q\left(\frac{n!}{(n-k)!}\right) = \frac{k+a(n-k)-a(n)}{q-1} \leq k + (q-1)\log n.
\end{displaymath}
Hence if $v_q(f_\lambda)\geq m+(q-1)\log n$, then $v_q(\chi^\lambda_\mu)\geq (m-k)$.
Thus taking $m=k+b$ in Theorem~A tells us that
\begin{displaymath}
  \lim_{n\to \infty} \frac{\#\{\lambda\vdash n \mid  v_p(\chi^\lambda_\mu)\leq b\}}{p(n)} = 0.
\end{displaymath}
From this the main theorem follows.
\bibliographystyle{abbrv}
\bibliography{refs}

\begin{thebibliography}{1}

\bibitem{MR2444214}
J.~Anderson.
\newblock An asymptotic formula for the {$t$}-core partition function and a
  conjecture of {S}tanton.
\newblock {\em J. Number Theory}, 128(9):2591--2615, 2008.

\bibitem{MR1321575}
A.~Granville and K.~Ono.
\newblock Defect zero {$p$}-blocks for finite simple groups.
\newblock {\em Trans. Amer. Math. Soc.}, 348(1):331--347, 1996.

\bibitem{HR}
G.~H. Hardy and S.~Ramanujan.
\newblock Asymptotic formulae in combinatory analysis.
\newblock {\em Proc. London Math. Soc. (2)}, 17(1):75--115, 1918.

\bibitem{MR2368985}
M.~Lassalle.
\newblock An explicit formula for the characters of the symmetric group.
\newblock {\em Math. Ann.}, 340(2):383--405, 2008.

\bibitem{Macdonald}
I.~G. Macdonald.
\newblock On the degrees of the irreducible representations of symmetric
  groups.
\newblock {\em Bull. London Math. Soc.}, 3:189--192, 1971.

\bibitem{Olsson}
J.~B. Olsson.
\newblock Mc{K}ay numbers and heights of characters.
\newblock {\em Math. Scand.}, 38(1):25--42, 1976.

\end{thebibliography}
\end{document}